\documentclass[12pt]{amsart}
\usepackage{amscd,times,fullpage}

\usepackage{amsthm,amsmath}
\usepackage{amssymb}
\usepackage{graphicx}
\usepackage{amsthm}
\usepackage{enumerate}
\usepackage{tikz-cd}
\usetikzlibrary{cd}
\usetikzlibrary{positioning}
\usetikzlibrary{matrix}

\newcommand{\Z}{\mathbb{Z}}

\newcommand{\CP}{\mathbb{C}\mathrm{P}}

\newcommand{\C}{\mathbb{C}}

\newcommand{\K}{K\"{a}hler }

\newcommand{\F}{\mathcal F}

\newcommand{\di}{{\operatorname{d}}}
\newcommand{\betti}{{b}}

\newcommand{\Id}{\operatorname{Id}}

\newcommand{\Sas}{\mathcal{S}}
\newcommand{\D}{\mathcal{D}}
\newcommand{\TM}{\text{T} M}

\newtheorem{theor}{Theorem}
\newtheorem{prop}[theor]{Proposition}

\newtheorem{lem}[theor]{Lemma}

\newtheorem{ex}[theor]{Example}
\newtheorem{remark}[theor]{Remark}

\makeatletter
\newtheorem*{rep@theorem}{\rep@title}
\newcommand{\newreptheor}[2]{%
\newenvironment{rep#1}[1]{%
 \def\rep@title{#2 \ref{##1}}%
 \begin{rep@theorem}}%
 {\end{rep@theorem}}}
\makeatother

\newreptheor{theor}{Theorem}

\usepackage{nicefrac}

\title{Sasaki structures distinguished by their basic Hodge numbers}

\author{D.~Kotschick}
\address{Mathematisches Institut, Ludwig-Maximilians-Universit\"at
M\"unchen, Theresienstr.~39, 80333 M\"unchen, Germany}
\email{dieter@math.lmu.de}

\author{G.~Placini}
\address{Dipartimento di Matematica e Informatica, Universit\`a di Cagliari, Via~Ospedale~72, 09124 Cagliari, Italy}
\email{giovanni.placini@unica.it}

\date{February 9, 2022 ; {\copyright \ D.~Kotschick and G.~Placini 2020}}

\subjclass[2010]{primary 53C25, 57R17; secondary 14M10, 14F45, 58A14}
\keywords{Sasaki structures, Hodge numbers, contact structures, complete intersections}

\begin{document}

\begin{abstract}
In all odd dimensions $\geq 5$ we produce examples of manifolds admitting pairs of Sasaki
structures with different basic Hodge numbers.  In dimension $5$ we prove more precise 
results, for example we show that on connected sums of copies of $S^2\times S^3$ the 
number of Sasaki structures with different basic Hodge numbers within a fixed homotopy 
class of almost contact structures is unbounded.
All the Sasaki  structures we consider are negative in the sense that the 
basic first Chern class is represented by a negative definite form of type $(1,1)$.
We also discuss the relation of these results to contact topology.
\end{abstract}
 
\maketitle

\section{Introduction}

A Sasaki structure consists of a contact form together with a compatible Riemannian metric and a 
transverse complex structure, which define a transverse K\"ahler structure for the Reeb foliation.
These structures are the odd-dimensional analogs of K\"ahler or projective-algebraic complex structures
on even-dimensional manifolds. 
The purpose of this paper is to exhibit families of closed manifolds which support two or more Sasaki
structures distinguished by their transverse Hodge numbers. Our most interesting examples are in dimension 
five, but we also prove results in arbitrary higher dimensions. 

There are several motivations for looking for manifolds with multiple Sasaki structures, some of them coming
from contact topology, and others coming from transverse K\"ahler geometry, involving not only the 
Hodge numbers but also the type or sign of the basic first Chern class. Before reviewing
these motivations, we briefly discuss the topological (non-)invariance of Hodge numbers on K\"ahler 
manifolds, particularly on complete intersections. This is of independent interest, and will be useful
for building examples in the Sasaki context.

\subsection{Hodge numbers of K\"ahler manifolds}

For compact K\"ahler manifolds the Hodge decomposition of the cohomology shows that certain linear 
combinations of Hodge numbers are topological invariants. The question whether there are additional 
linear combinations that are topological invariants was asked first by Hirzebruch in his famous collection of open problems 
from 1954, but was only resolved less than ten years ago by the first named author and Schreieder~\cite{KS}.
Note that for questions like these the intuition gained from small dimensions may be misleading, since in 
complex dimension $1$ all Hodge numbers are determined by Betti numbers, and in complex dimension $2$ 
they are determined by the Betti numbers and the signature. However, starting in dimension $3$, there are 
always more Hodge numbers than there are obvious topological invariants. Indeed the answer to Hirzebruch's
question obtained in~\cite{KS} is that the only linear combinations of Hodge numbers which are invariant 
under orientation-preserving diffeomorphisms between compact K\"ahler manifolds are those that can be 
expressed through the Betti numbers and the signature.  
 
After the results of~\cite{KS} the question arose whether more Hodge numbers or combinations thereof might 
be invariant if one restricts the class of K\"ahler manifolds one considers, for example, by looking only at 
manifolds with a simple enough cohomology structure. In this direction, Chataur~\cite{Chat} asked whether 
diffeomorphic complete intersections must always have the same Hodge numbers. The first dimension for which 
this question is interesting is dimension $3$. It is not hard to see from the Hirzebruch--Riemann--Roch theorem 
that two three-dimensional complete intersections have different Hodge numbers if and only if they have 
different first Chern classes; see Lemma~\ref{LemHodgeChernCompInt} below. 
Since examples with distinct first Chern classes were found by Wang and Du~\cite{wangdu16}, 
this answers Chataur's question in the negative. After this work was completed, we found that this conclusion, 
based on the same examples, was also reached recently by Wang, Yu and Wang~\cite{wangetal}.

\subsection{Basic Hodge numbers of Sasaki manifolds}

For a Sasaki structure the orbits of the Reeb flow are the leaves of a Riemannian foliation, for which
one can consider the basic cohomology, giving rise to the notion of basic Betti numbers. 
The properties of the basic cohomology were studied by El~Kacimi-Alaoui, Hector and 
Nicolau~\cite{elkacimialaoui90,elkacimialaouihector86,elkacimialaouinicolau93}. One of their 
conclusions is that the basic Betti numbers are finite, and are topological invariants of the underlying 
manifold, so they do not depend on the choice of Sasaki structure; 
compare~\cite[Theorem~7.4.14]{boyergalicki08}.

Since the Reeb foliation is transversely K\"ahler, the basic cohomology has a Hodge decomposition
mimicking precisely what happens for K\"ahler manifolds. Once again there are linear combinations 
of Hodge numbers which can be expressed through the Betti numbers, and are therefore topological
invariants. One can now ask for these basic Hodge numbers of Sasaki manifolds 
the analog of Hirzebruch's question for compact K\"ahler manifolds which we discussed above.

Concerning this question, it has been known for some time that the basic Hodge numbers are not always 
topological invariants. However, up to now, there has been only a single example of two Sasaki structures 
on the same manifold but with different basic Hodge numbers. This example was alluded to by Boyer and 
Galicki in their monograph~\cite[p.~234]{boyergalicki08}, but did not appear in print until the work of  
Goertsches, Nozawa and T\"oben~\cite{goertschesnozawatoeben16}. In that example, 
see~\cite[Example~3.4]{goertschesnozawatoeben16}, one 
of the Sasaki structures is positive, in the sense that its basic first Chern class is represented 
by a positive $(1,1)$-form, which implies the vanishing of certain Hodge numbers. The other 
Sasaki structure is null rather than positive, and has different Hodge numbers. This raises the
further question whether different Hodge numbers can only occur for pairs of Sasaki structures 
with different type.

Very recently, Ra\'{z}ny~\cite{razny19} proved that the basic Hodge numbers of Sasaki
manifolds are constant in arbitrary smooth deformations of the Sasaki structure. This completes
earlier results of Boyer and Galicki~\cite{boyergalicki08} and of Goertsches, Nozawa and 
T\"oben~\cite{goertschesnozawatoeben16} who proved the invariance for special types
of deformations. Ra\'{z}ny's result gives further impetus to the study of basic Hodge numbers,
since it shows they provide a tool, independent of contact topology, to distinguish deformation
classes of Sasaki structures.

We show in this paper that in every odd dimension $\geq 5$ there are manifolds which have Sasaki 
structures with different Hodge numbers. Such examples exist in infinitely many distinct 
homotopy types in each dimension. Moreover, unlike in the example 
from~\cite{boyergalicki08,goertschesnozawatoeben16} mentioned above, they can be chosen to 
be of the same type, so that the difference in Hodge numbers is not explained by the type or sign 
of the first Chern class. We use Sasaki structures of negative type, which are of course plentiful
as they can be constructed from complex projective manifolds with ample canonical bundle.

Given the results of~\cite{KS} discussed above, it should not be surprising that the basic Hodge numbers 
of diffeomorphic Sasaki manifolds may disagree. Indeed, whenever we have two diffeomorphic closed 
K\"ahler manifolds with integral K\"ahler classes, we can try to build diffeomorphic 
Sasaki manifolds by the so-called Boothby--Wang construction~\cite{boothbywang58}, using circle 
bundles whose Euler classes are the K\"ahler classes. For these special Sasaki structures the basic 
Hodge numbers agree with the usual Hodge numbers of the base manifolds, so if the two K\"ahler 
manifolds one starts with have different Hodge numbers, then so do the Sasaki structures on the 
circle bundles. The difficulty in turning this heuristic into a rigorous argument is that the diffeomorphism of K\"ahler 
manifolds one starts with will only lift to a diffeomorphism of the total spaces of the circle bundles
if it maps one Euler class to the other, that is, it sends one K\"ahler class to the other. For manifolds 
with large second Betti number this is a condition which can be hard to check, since one has to 
compare the K\"ahler cones of distinct complex structures which are completely unrelated to each other.
However, if the K\"ahler manifolds have $b_2=1$, then this difficulty does not arise. Therefore the 
examples of complex three-dimensional diffeomorphic complete intersections with distinct Hodge numbers 
mentioned above do indeed give examples of Sasaki structures on diffeomorophic $7$-manifolds but with 
distinct Hodge numbers.

This kind of argument does not produce $5$-dimensional examples, because for complex surfaces 
the Hodge numbers are diffeomorphism invariants, cf.~\cite[Theorem~2]{BLMSorient}. However, the classification of simply 
connected $5$-manifolds due to Smale~\cite{smale} and Barden~\cite{barden65} shows that we 
can build diffeomorphic $5$-manifolds as circle bundles over $4$-manifolds which need not be 
diffeomorphic. In particular, we will see that we can take the $4$-manifolds to be algebraic surfaces 
with different Hodge numbers. 

\subsection{Symplectic topology}

Eliashberg \cite{eliashberg89} introduced a dichotomy of $3$-dimensional contact structures into overtwisted 
and tight ones, and proved that any homotopy class of almost contact structures contains a unique isotopy 
class of overtwisted contact structures. After many attempts to extend the notion of overtwistedness to higher 
dimensions, a satisfactory definition was finally given by Borman, Eliashberg and 
Murphy~\cite{bormaneliashbergmurphy15}, who also proved the classification via an $h$-principle in all dimensions.

The contrast between overtwisted and tight structures is encountered, for instance, when considering symplectic 
fillability. Both in dimension three and higher the existence of a symplectic filling gives an obstruction to 
overtwistedness~\cite{eliashberg90,niederkrueger06}.
As a consequence of a result of Niederkr\"uger and Pasquotto \cite{niederkruegerpasquotto09}  the contact structures
underlying Sasaki structures are symplectically fillable, hence tight. Therefore they do not fall into the classification given 
by an $h$-principle, so that it may happen that two Sasaki structures on the same manifold can have homotopic underlying 
almost contact structures but non-isotopic contact structures. Examples of this phenomenon have been exhibited by Boyer, 
Macarini and van Koert~\cite{boyermacarini}, who showed that it occurs for positive Sasaki structures on 
connected sums of copies of $S^2\times S^3$.

Note that the deformations of Sasaki structures considered by Ra\'{z}ny~\cite{razny19} preserve the isotopy
class of the underlying contact structure. It is not clear how the basic Hodge numbers of a Sasaki 
structure relate to the underlying contact structure; cf.~the discussion in Remark~\ref{remm}.

\subsection{Statement of results}

The first dimension in which the above questions are interesting for Sasaki manifolds is dimension $5$, since in
dimension $3$ the basic Hodge numbers are determined by the Betti numbers. 
For simply connected $5$-manifolds we have the classification due to Smale~\cite{smale} and Barden~\cite{barden65},
which in particular allows us to find the same $5$-manifold as the total space of Boothby--Wang fibrations over 
very different algebraic surfaces. In this way many interesting examples can be constructed. Instead of 
stating comprehensive results for many different manifolds, we focus on the simplest case $M_n=n(S^2\times S^3)$, 
the connected sum of $n$ copies of $S^2\times S^3$.
For these manifolds we prove the following result which exhibits the topological non-invariance of the basic 
Hodge numbers:
\begin{theor}\label{main1}
\begin{itemize}
\item[(A)] For every $k$ there is an $n$ such that $M_n$ has at least $k$ negative Sasaki structures whose 
underlying almost contact structures are homotopic, but whose contact structures are pairwise inequivalent.
\item[(B)] For every $k$ there is an $n$ such that $M_n$ has at least $k$ negative Sasaki structures whose 
underlying almost contact structures are homotopic, but which have pairwise distinct Hodge numbers.
\item[(C)] For every $n=80m+76$, the manifold $M_n$ has two negative Sasaki structures whose 
underlying almost contact structures are homotopic, but whose contact structures are inequivalent,
and whose Hodge numbers disagree.
\end{itemize}
\end{theor}
Negativity means that the basic first Chern class is represented by a negative  $(1,1)$-form.
It is a result of Gomez~\cite{gomez} that $M_n$ admits a negative Sasaki structure for every $n$.
By construction, the first Chern class of the underlying contact structure vanishes for all our examples.
This means that the underlying almost contact structure is the same as the one for the positive 
Sasaki structures considered by Boyer, Macarini and van Koert~\cite{boyermacarini}.

A few more remarks on this theorem are in order. While part (A) is a straightforward application of 
known results, it does not give examples with distinct basic Hodge numbers. To prove part (B), involving the 
Hodge numbers, 
we actually construct a convenient supply of algebraic surfaces to be used as bases for the 
Boothby--Wang fibrations, and this construction is of independent interest. It is clear that in statement
(B) the number $n$ must grow as we increase $k$.  While it would be possible to tweak the proof of 
(B) so as to control particular invariants of the contact structure, it seems impossible to actually prove 
that there are such examples with equivalent contact structures, cf.~Remark~\ref{remm}. 
In both (A) and (C) the inequivalence of 
the contact structures arises from a result of Hamilton~\cite{hamilton13}, who implemented an idea
of Eliashberg, Givental and Hofer~\cite{EGH}.

The phenomena exhibited in Theorem~\ref{main1} also occur on non-spin manifolds. The simplest 
examples of these are $N_n=M_n\# (S^2\widetilde\times S^3)$, where $S^2\widetilde\times S^3$
denotes the non-trivial orientable $S^3$-bundle over $S^2$. Concerning the basic Hodge numbers we 
will prove the following result:
\begin{theor}\label{5nonspin}
For every $k$ there is an $n$ such that $N_n$ has at least $k$ negative Sasaki structures with
pairwise distinct Hodge numbers.
\end{theor}

Finally we show that the topological non-invariance of basic Hodge numbers is true in higher 
dimensions as well. 
\begin{theor}\label{main2}
For all $n>1$ there exist $(2n+1)$-dimensional manifolds which have pairs of negative Sasaki 
structures with different basic Hodge numbers. Moreover, in each dimension manifolds with 
such pairs range over infinitely many homotopy types.
\end{theor}
Except for the cases $n=3$ and $4$, the proof will show that all these manifolds may be taken to be 
simply connected.

We have taken care to arrange all the examples in the proofs of these theorems to have negative 
Sasaki structures. Some of the constructions would be more flexible if one did not pay attention
to the sign of the first Chern class, and in particular if one were willing to use indefinite Sasaki 
structures.

\section{Sasaki structures and their transverse geometry}\label{sectiontop}

We begin with some definitions and known results, for a more exhaustive treatment we refer to the monograph 
by Boyer and Galicki~\cite{boyergalicki08}. 
All manifolds are assumed to be smooth, connected, closed and oriented.

A \textit{contact form} $\eta$ on a manifold $M$ of dimension $2n+1$ is a $1$-form with the property that 
$\eta\wedge (d\eta)^n$ is a volume form, positive with respect to the given orientation. This means 
that the \textit{contact structure} $\ker (\eta )$ is a hyperplane field endowed with the positive
symplectic form $d\eta$. Therefore, a contact form defines a reduction of the structure group
of the tangent bundle from $SO(2n+1)$ to $\{ 1 \}\times U(n)$, and we refer to (the homotopy class of)
this reduction as the underlying \textit{almost contact structure}. 

A \textit{K-contact structure} $(\eta,\phi,R,g)$ on $M$ consists of a contact form $\eta$ and an endomorphism 
$\phi$ of the tangent bundle $\TM$ satisfying the following properties:
\begin{enumerate}
\item[$\bullet$] $\phi^2=-\Id+R\otimes\eta$ where $R$ is the Reeb vector field of $\eta$,
\item[$\bullet$] $\phi_{\vert\D}$ is an almost complex structure compatible with the symplectic form $\di\eta$ on $\D=\ker\eta$,
\item[$\bullet$] the Reeb vector field $R$ is Killing with respect to the metric $g(\cdot,\cdot)=\di\eta(\phi\cdot,\cdot)+\eta(\cdot)\eta(\cdot)$.
\end{enumerate}   
Given such a structure one can consider the almost complex structure $I$ on the Riemannian cone $\big( M\times(0,\infty),t^2g+\di t^2\big)$ given by
\begin{enumerate}
\item[$\bullet$] $I=\phi$ on $\D=\ker\eta$, and
\item[$\bullet$] $I(R)=t\partial_t$.
\end{enumerate} 
A \textit{Sasaki structure} is a K-contact structure $(\eta,\phi,R,g)$ such that the associated almost complex structure $I$ is integrable. 
We call a manifold $M$ \textit{Sasakian} if it admits a Sasaki structure.
A \textit{Sasaki} manifold is a manifold equipped with a Sasaki structure.

A Sasaki manifold is called \textit{regular} (respectively \textit{quasi-regular}, \textit{irregular}) if its Reeb foliation has the corresponding
property.
Every regular Sasaki manifold is a \textit{Boothby--Wang fibration} $M$ over a projective manifold $(X,\omega)$ with $\omega$ representing an integral class (\cite{boothbywang58,boyergalicki08}), that is, the principal $S^1$-bundle $\pi\colon M\longrightarrow X$ with Euler class $[\omega]$ and connection $1$-form $\eta$ such that $\pi^*(\omega)=\di\eta$.
All the Sasaki structures we consider in our examples will be of this type.

The Reeb foliation $\F$ of a Sasaki structure is transversally K\"ahler. 
Cohomological properties of the transverse space, that is, the so-called basic cohomology of the foliation, 
were studied by El~Kacimi-Alaoui, Hector and Nicolau~\cite{elkacimialaoui90,elkacimialaouihector86,elkacimialaouinicolau93}.
They proved many transverse analogues of classical properties of \K manifolds.
Namely, one can define a basic Dolbeault double complex  $\Omega^{\bullet,\bullet}_B(\F)$ and prove that it satisfies 
the Hodge decomposition theorem and Poincar\'e and Serre dualities. 
Moreover, in parallel with the standard case one can define basic Chern classes $c_i(\F)$ with the use of a transverse connection.
It is then natural to call a Sasaki structure \textit{positive}, respectively \textit{negative} or \textit{null}, if $c_1(\F)$ can be represented 
by a positive definite, resp.~negative definite or null,  basic $(1,1)$-form. It is called indefinite otherwise.
Again in analogy with the K\"ahler case one defines the \textit{basic Betti and Hodge numbers} $\betti_B^r(\F)$ and $h^{p,q}_B(\F)$
 to be the dimensions of basic de Rham and Dolbeault cohomology groups respectively.

If $M$ is a regular Sasaki manifold, basic forms on $M$ corresponds exactly to forms on the base $X$ of the Boothby--Wang 
fibration $\pi\colon M\longrightarrow X$. Moreover, the basic Dolbeault complex $\Omega^{\bullet,\bullet}_B(\F)$ is the usual 
Dolbeault complex $\Omega^{\bullet,\bullet}(X)$ of $X$. Thus, the basic Betti and Hodge numbers $\betti_B^r(\F)$ and 
$h^{p,q}_B(\F)$ are the Betti and Hodge numbers $\betti^r(X)$ and $h^{p,q}(X)$ of $X$.
Also, the basic Chern classes $c_i(\F)$ are the Chern classes $c_i(X)$ of $X$.
Therefore, the Sasaki structure on $M$ is positive, respectively negative or null, if and only if the first Chern class $c_1(X)$ can 
be represented by a positive definite, resp.~negative definite or null, form of type $(1,1)$.

The following result shows that basic Betti numbers are topological invariants of Sasakian manifolds.
\begin{theor}[{\cite[Theorem~7.4.14]{boyergalicki08}}]\label{ThmInvBasicBettiNumbers}
Let $(M,\eta,\phi,R,g)$ be a compact Sasaki manifold. The basic cohomology $H^*_B(\F)$ only depends on the topology of $M$. 
In particular, the basic Betti numbers of any two Sasaki structures on $M$ agree. 
\end{theor}
There can be no such result for the basic Hodge numbers, which are geometric invariants and can distinguish Sasaki structures. 
The first example of this was given in~\cite[Example~3.4]{goertschesnozawatoeben16}.

\begin{ex}\label{ExGoertsches}
Consider $M_{21}$, the connected sum of $21$ copies of $S^2\times S^3$. 
On the one hand, this manifold can be endowed with the null Sasaki structure $\Sas_1$ given by 
the Boothby--Wang fibration over a $K3$ surface. 
Since this structure is regular, the basic Hodge numbers of $\Sas_1$ are the Hodge numbers of the $K3$ surface. 
In particular $h^{2,0}(\F_1)=1$.
On the other hand, $M$ supports a positive Sasaki structure $\Sas_2$ arising as follows; see \cite[page~356]{boyergalicki08}. 
The connected sum $\#21(S^2\times S^3)$ can be realized as the link $L_f=V_f\bigcap S^7\subset \C^4$ where 
\begin{align*}
f(z_0,z_1,z_2,z_3)= z_0^{22}+z_1^{22}+z_2^{22}+z_0z_3
\end{align*}
is a weighted homogeneous polynomial of degree $d=22$ with weight $w=(1,1,1,21)$.
This Sasaki structure is positive because $\sum w_i-d=2>0$, see \cite[Proposition~9.2.4]{boyergalicki08}.
Then a vanishing result proved independently by Nozawa~\cite{nozawa14} and by Goto~\cite{goto12}  guarantees that $h^{2,0}(\F_2)=0$. 
\end{ex}

\section{Some complex projective manifolds with distinct Hodge numbers}\label{sectionhypersurfaces}

Since for algebraic surfaces the Hodge numbers are diffeomorphism invariants~\cite{BLMSorient}, there is 
no hope for proving Theorem~\ref{main1} (B) using the Boothby--Wang construction on diffeomorphic 
algebraic surfaces. However, in order to obtain diffeomorphic
$5$-manifolds via the Boothby--Wang construction it is enough to ask that the surfaces we start with
are simply connected and have the same Euler characteristic. A convenient supply of examples is 
provided by the following theorem.

\begin{theor}\label{surfaces}
For every positive integer $k$ there are $k$-tuples 
%$X_1, \ldots , X_k$ 
of simply connected algebraic surfaces 
with ample canonical bundles which have the same Euler characteristic, but have pairwise distinct Hodge 
numbers.
\end{theor}
In terms of Chern numbers, the conditions are that all the surfaces in such a tuple have the same $c_2$,
but have pairwise distinct $c_1^2$, equivalently, pairwise distinct signatures.
\begin{proof}
Our examples will be generic smooth hypersurfaces of bidegree $(p,3q)$ in $\C P^1\times \C P^2$. For $p>2$ and $q>1$
they have ample canonical bundles. Moreover, they are simply connected by the Lefschetz hyperplane theorem.
By the adjunction formula the first Chern class of such a hypersurface is 
\begin{equation}\label{c1}
c_1 = (2-p)x_1+(3-3q)x_2 \ ,
\end{equation}
where the $x_i$ are the generators of the second cohomology coming from the two factors of the product $\C P^1\times \C P^2$.
A straightforward calculation yields the following Chern numbers:
\begin{equation}\label{c12}
c_1^2 = 9(q-1)(3pq-p-4q) \ ,
\end{equation}
\begin{equation}\label{c2}
c_2 = 3\big( p(3q-1)^2-6q(q-1)\big) \ .
\end{equation}
The sequence of integers $3q-1$ contains infinitely many primes, so {\it a fortiori} it is possible to choose 
arbitrarily large sets of such numbers which are pairwise coprime. Let us start with a set of $3k$ values of $q$
for which the numbers $3q-1$ are pairwise coprime.
For such a set we want to find an integer $n$
so that $c_2=3n$ will be realised by all our choices of $q$. This means that $n$ must satisfy
\begin{equation}\label{cong}
n \equiv 6q(1-q) \mod{(3q-1)^2}
\end{equation}
for every $q$. Since we chose the numbers $3q-1$ to be pairwise coprime, the Chinese Remainder Theorem
guarantees the existence of solutions to this system of congruences for all our $q$ simultaneously. 

For an $n$ solving all the congruences~\eqref{cong}, we can solve the equation~\eqref{c2} with 
$c_2=3n$ for every $q$ to obtain a unique positive integer $p$ such that the hypersurface 
of bidegree $(p,3q)$ in $\C P^1\times \C P^2$ has Euler number $3n$. This $p$ is given by 
$$
p=\frac{n+6q(q-1)}{(3q-1)^2} \ ,
$$
and substituting this into~\eqref{c12} we obtain
$$
c_1^2 = \frac{9(q-1)}{3q-1}\big( n-2q(3q+1)\big) \ .
$$
For a fixed value of $c_1^2$ this becomes a cubic equation for $q$, showing that among our $3k$
surfaces at least $k$ distinct values for $c_1^2$ are realised. This finishes the proof of the theorem.
\end{proof}

\begin{ex}\label{five}
To illustrate the theorem, let us start with $q\in\{ 2,3,4,6,8\}$. Then the $3q-1$ range over $5, 8, 11, 17, 23$
and so are indeed pairwise coprime. The smallest positive $n$ solving the system of five congruences is 
$$
n= 21.740.924.188 \ ,
$$
from which one computes the values of $p$ and $c_1^2$ for every $q$. The results are shown in 
Table~\ref{TableSurfaces}.
\begin{table}[ht]
\centering
\begin{tabular}{c| c| c| c} 
$q$ & $p$ & $c_1^2$ &  $d (c_1)$ \\ [0.5ex] 
\hline %\rule{0pt}{3ex} 
$2$ & $869.636.968$ & $39.133.663.488$ &  $1$ \\ % inserting body of the table
$3$ & $339.701.941$ & $48.917.079.288$ &  $1$ \\
$4$ & $179.677.060$ & $53.364.086.388$ &  $1$ \\
$6$ & $75.228.112$ & $57.549.504.600$ &  $5$ \\
$8$ & $41.098.156$ & $59.551.226.028$ &  $1$ \\
\hline
\end{tabular}
\vspace{10pt}
\caption{A five-tuple of surfaces as in Theorem~\ref{surfaces}} % title of Table
\label{TableSurfaces}
\end{table}
Here there are no coincidences among the values of $c_1^2$, so we obtain five surfaces with the same $c_2$ and 
with pairwise distinct Hodge numbers. The final column in the table shows the divisibility of the first Chern class.
\end{ex}

In higher dimensions the Hodge numbers are not always diffeomorphism-invariant~\cite{KS},
so other constructions are possible.

Let $X$ be a complete intersection of complex dimension $3$. As a consequence of the Lefschetz hyperplane theorem 
$X$ is simply connected and $h^{p,q}(X)=h^{p,q}(\CP^{3})$ for $p+q<3$. Moreover, as observed by Hirzebruch a long 
time ago, the Hodge numbers are actually determined by the Chern numbers in this case; compare~\cite[Appendix One]{hirzebruch95}. 
The following includes a strong converse to this observation.

\begin{lem}\label{LemHodgeChernCompInt}
Two diffeomorphic complete intersections $X$ and $Y$ of complex dimension $3$ have different Hodge numbers if and only if $c_1(X)\neq c_1(Y)$.
\end{lem}
\begin{proof}
In this case the only interesting Hodge numbers are $h^{0,3}=h^{3,0}$ and $h^{1,2}=h^{2,1}$. They determine each other by
the equation 
\begin{equation}\label{eq:Hodge}
2h^{0,3}+2h^{1,2}=b_3=4-c_3 \ .
\end{equation}
By the Hirzebruch--Riemann--Roch Theorem~\cite{hirzebruch95}  we have:
\begin{equation}\label{eq:HRR}
h^{0,3}=1-\chi=1-\frac{1}{24}c_1c_2 \ ,
\end{equation}
where $\chi=\chi (\mathcal{O})$ is the holomorphic Euler characteristic.

Let $X$ be a smooth complete intersection of codimension $r$ with multi-degree $(d_1,\ldots, d_r)$. Without loss 
of generality we assume $d_i\geq 2$ for all $i$. We have $c_1(X)=k\cdot x$ and $p_1(X)=m\cdot x^2$, where 
$x$ is the positive generator of the second cohomology, and 
$$
k=4+r-\sum_{i=1}^r d_i
$$
and
$$
m=4+r-\sum_{i=1}^r d_i^2 \ .
$$
The last equation shows that $m$ is negative, unless $r=1$ and $X$ is a quadric in $\C P^4$.
So there can be no other complete intersection diffeomorphic to a quadric.

If $Y$ is another complete intersection that is diffeomorphic to $X$, then, up to conjugating the complex 
structure, we may assume that the diffeomorphism pulls back the positive generator in cohomology to the 
positive generator. Therefore, $X$ and $Y$ have the same first Chern class if and only if they have 
the same $k$. Moreover, they have the same $m$ just because they are diffeomorphic.
For a fixed $k$ we have 
\begin{align*}
2c_2&=c_1^2-p_1=(k^2-m)\cdot x^2
\end{align*}
and 
\begin{align*}
2c_1c_2&=k(k^2-m)\cdot x^3\ .
\end{align*}
Together with equations~\eqref{eq:Hodge} and~\eqref{eq:HRR}, this shows that if $X$ and $Y$ are 
diffeomorphic and have the same first Chern class, then their Hodge numbers agree.

Conversely, suppose that $X$ and $Y$ are diffeomorphic, and so have the same $m$, and also 
have the same Hodge numbers. Then by~\eqref{eq:HRR} they have the same $2c_1c_2$ which
is given by the formula 
$$
f(k) = k(k^2-m)
$$
as a function of the first Chern class $k\cdot x$. The derivative of this function is $3k^2-m$.
Since by our earlier observation about quadrics we may assume that $m$ is negative, the 
derivative of $f$ is strictly positive, showing that the function is strictly monotonically increasing.
In other words, if $X$ and $Y$ have the same Hodge numbers, then they have the same first Chern classes.
\end{proof}

By Wall's classification the diffeomorphism type of a $3$-dimensional complete intersection is determined by the 
Euler number, the first Pontryagin class and the parity of the first Chern class, see~\cite{jupp73,wall66}.
Libgober and Wood~\cite{libgoberwood82} conjectured the existence of diffeomorphic complete intersections 
with different Chern classes. However, the first such pairs were only found comparatively recently by Wang and Du~\cite{wangdu16},
we list them in Table~\ref{TableDiffeoHodgeNumbers}. By Lemma~\ref{LemHodgeChernCompInt} above, these pairs have distinct 
Hodge numbers. 
\begin{table}[ht]
\centering
\begin{tabular}{c| c c c c} 
$\underline{d}$ & $d$ & $p_1$ & $\chi=d\cdot c_3$ & $c_1$ \\ [0.5ex] 
\hline \rule{0pt}{3ex} 
$(70,16,16,14,7,6)$& $7^3\cdot5\cdot3\cdot2^{11}$ & $-5683$ & $-7767425433600$ & $-119$ \\ % inserting body of the table
$(56,49,8,6,5,4,4)$& $7^3\cdot5\cdot3\cdot2^{11}$ & $-5683$ & $-7767425433600$ & $-121$ \\\hline \rule{0pt}{3ex} 
$(88,28,19,14,6,6)$ & $19\cdot11\cdot7^2\cdot3^2\cdot2^8$ & $-9147$ & $-35445749391360$ & $-151$ \\
$(76,56,11,7,6,6,2)$ & $19\cdot11\cdot7^2\cdot3^2\cdot2^8$ & $-9147$ & $-35445749391360$ & $-153$ \\\hline \rule{0pt}{3ex} 
$(84,29,25,25,18,7)$ & $29\cdot7^2\cdot5^4\cdot3^3\cdot2^3$ & $-9510$ & $-384536710530000$ & $-178$ \\
$(60,58,49,9,5,5,5)$ & $29\cdot7^2\cdot5^4\cdot3^3\cdot2^3$ & $-9510$ & $-384536710530000$ & $-180$  \\ [1ex] % [1ex] adds vertical space
\hline %inserts single line
\end{tabular}
\vspace{10pt}
\caption{Diffeomorphic $3$-dimensional complete intersections with different $c_1$} % title of Table
\label{TableDiffeoHodgeNumbers}
\end{table}
The characteristic classes of a complete intersection $X$ are multiples of the generators $x,x^2$ of the groups $H^2(X)\cong H^4(X)\cong\Z$.
The values in Table~\ref{TableDiffeoHodgeNumbers} are the coefficients that determine the characteristic classes as multiples of $x$ and $x^2$.

\begin{remark}
The conclusion that the diffeomorphic complete intersection threefolds with different first Chern classes
from~\cite{wangdu16} have distinct Hodge numbers was also reached recently in~\cite{wangetal}. There, the 
Hodge numbers are computed by brute force, without a general result along the lines of Lemma~\ref{LemHodgeChernCompInt}.
The paper~\cite{wangetal} also contains a pair of five-dimensional diffeomorphic complete intersections 
with distinct Hodge numbers.
\end{remark}

\section{Five-dimensional Sasaki manifolds}\label{sectionproofs0}

In this section we discuss five-dimensional Sasaki manifolds obtained as Boothby--Wang fibrations over 
algebraic surfaces. In particular, we prove Theorems~\ref{main1} and \ref{5nonspin}.

Let $X$ be a smooth simply connected algebraic surface, and $[\omega]\in H^2(X,\Z)$ an integral K\"ahler class.
We consider the Sasaki manifold $M$ obtained as the total space of the Boothby--Wang fibration $\pi\colon M\longrightarrow X$
with Euler class $[\omega]$. Then the first Chern class of the almost contact structure underlying the 
Sasaki structure is $\pi^*(c_1(X))$.

Recall that the \textit{divisibility} $d(\alpha)$ of a class $\alpha\in H^2(X;\Z)$ is the maximum number $n\in\Z$ such that 
$\alpha=n\beta$ for some $0\neq\beta\in H^2(X;\Z)$. A class $\alpha$ is called \textit{indivisible}, or \textit{primitive}, if $d(\alpha)=1$.
For our construction we always scale the K\"ahler form $\omega$ so that it represents a primitive cohomology
class. We can then use the following standard result, proved for example in the paper of Hamilton~\cite{hamilton13}.
\begin{lem}\label{lem:spin}
If $X$ is simply connected and $[\omega]$ is primitive, then the total space $M$ of the Boothby--Wang 
fibration $\pi\colon M\longrightarrow X$ is simply connected with torsion-free cohomology. It is spin if
and only if $X$ is spin or $c_1(X)\equiv [\omega]\mod 2$.
\end{lem}
\begin{proof}
The statements about the fundamental group and the cohomology follow from the exact homotopy sequence and
the Gysin sequence for $\pi$, see~\cite[Section~4]{hamilton13} for the details. Moreover, $w_2 (M)$ is the mod $2$ 
reduction of $\pi^*(c_1(X))$. This vanishes if $X$ is spin, since then $c_1(X)$ has trivial mod $2$ reduction, or, 
more generally, if and only if the mod $2$ reduction of $c_1(X)$ is in the kernel of $\pi^*$, which is spanned by the 
Euler class $[\omega]$; compare~\cite[Lemma~26]{hamilton13}.
\end{proof}
Combining this with the classification of simply connected $5$-manifolds and of almost contact structures
on them, we obtain:
\begin{prop}\label{SmaleGeiges}
Suppose that $c_1(X)$ is a (positive or negative) multiple of $[\omega]$.
Then the total space $M$ is diffeomorphic to $n(S^2\times S^3)$ with $n=b_2(X)-1$. The almost contact 
structure underlying the Sasaki structure is the unique one with trivial first Chern class.
\end{prop}
\begin{proof}
The assumption that $c_1(X)$ is a multiple of the Euler class $[\omega]$ implies that $M$ is spin.
The classification of simply connected spin $5$-manifolds due to Smale~\cite{smale} implies that 
$M$ is diffeomorphic to $n(S^2\times S^3)$ with $n=b_2(X)-1$. Moreover, $c_1(X)$ is in the 
kernel of $\pi^*$, and so the first Chern class of the (almost) contact structure vanishes. It is a result
of Geiges~\cite{geiges91} that in this case the first Chern class is a complete invariant of 
almost contact structures.
\end{proof}

Recall that two contact structures are called \textit{equivalent} if they can be related by a series of diffeomorphisms and isotopies.
Within the same homotopy class of almost contact structures we can sometimes distinguish inequivalent contact structures using
the following special case of a theorem of Hamilton~\cite[Corollary~43]{hamilton13}.
\begin{prop}[\cite{hamilton13}]\label{ThmHamilton}
Let $M$ be a simply connected $5$-manifold admitting two different regular Sasaki structures $(\eta_i,\phi_i,R_i,g_i)$ for $i=1,2$ with indivisible Euler classes $[\omega_i]$, viz.
$$
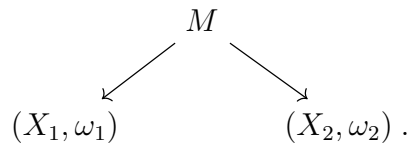
\begin{tikzcd}[column sep=small]
& M\arrow[dl]\arrow[dr] & \\
(X_1,\omega_1) && (X_2,\omega_2)\ .
\end{tikzcd}
$$
If the contact structures defined by $\eta_1$ and $\eta_2$ have trivial first Chern class and are equivalent, then $d(c_1(X_1))=d(c_1(X_2))$.
\end{prop}

We now prove Theorem~\ref{main1}.
\begin{proof}[Proof of Theorem~\ref{main1} (A)]
In~\cite{braungardtkotschick05} Braungardt and the first named author constructed arbitrarily large tuples 
$X_1,\ldots,X_k$ of homeomorphic branched covers of the projective plane $\CP^2$. 
These are pairwise non-diffeomorphic projective surfaces distinguished by the divisibilities of their first Chern classes. 
Moreover, they all have ample canonical bundles; cf.~\cite[Corollary~1]{braungardtkotschick05}. 

On each of the surfaces in such a $k$-tuple we perform the Boothby--Wang construction using as Euler 
classes the primitive integral cohomology classes which are rational multiples of the canonical classes.
By Proposition~\ref{SmaleGeiges} this yields $k$-tuples of negative Sasaki structures on $M_n=n(S^2\times S^3)$
with $n=b_2(X_i)-1$, which are homotopic as almost contact structures.

However, since the  first Chern classes $c_1(X_i)$ have pairwise different divisibilities, the
contact structures are inequivalent by Proposition~\ref{ThmHamilton}. 
\end{proof}

\begin{proof}[Proof of Theorem~\ref{main1} (B)]
We start with a $k$-tuple of surfaces $X_1,\ldots X_k$ as in Theorem~\ref{surfaces}, 
so they are simply connected with ample canonical bundles, with the same Euler characteristics,
but with pairwise distinct Hodge numbers.

Again we perform the Boothby--Wang construction on each surface using as Euler 
classes the primitive integral cohomology classes which are rational multiples of the canonical classes.
By Proposition~\ref{SmaleGeiges} this yields $k$-tuples of negative Sasaki structures on $M_n=n(S^2\times S^3)$
with $n=b_2(X_i)-1$, which are homotopic as almost contact structures. Their basic Hodge numbers
 are the usual Hodge numbers of the $X_i$, so for $i\neq j$ the Sasaki structures obtained from $X_i$
 and from $X_j$ have different Hodge numbers.
\end{proof}
\begin{remark}\label{remm}
In general we do not have any control over the contact structures in these examples. 
For instance, in the five-tuple considered in Example~\ref{five}, two different divisibilities for the 
first Chern classes occur, so by Proposition~\ref{ThmHamilton} there are at least two different contact 
structures underlying these Sasaki structures. If one wants to prove a version of Theorem~\ref{main1} (B)
where one has tuples of Sasaki structures with pairwise distinct Hodge numbers, but such that all the 
contact structures are equivalent, then one could try to sharpen our proof above, to ensure that all the 
divisibilities are the same, so that the contact structures could not be distinguished 
using Proposition~\ref{ThmHamilton}. However, the vanishing of this particular obstruction does 
not prove that the contact structures really are equivalent.
\end{remark}

\begin{proof}[Proof of Theorem~\ref{main1} (C)]
These structures are obtained as Boothby--Wang fibrations over a family of complete intersections and a family of 
Horikawa surfaces respectively.

On the one hand, let $X_k$ be a generic complete intersection in $\CP^1\times\CP^3$ given by intersecting 
hypersurfaces of  bidegree $(2,5)$ and $(k,1)$.  
By the adjunction formula we have $c_1(X_k)=-kx_1-2x_2$ where $x_1$ and $x_2$ are the generators of the 
cohomology rings of $\CP^1$ and $\CP^3$ respectively.
The Chern numbers and holomorphic Euler characteristic of $X_k$ are:
$$
c_1^2(X_k)=40k+8 \ ,\hspace{3mm}c_2(X_k)=80k+76 \ ,\hspace{3mm}\chi(\mathcal{O}_{X_k})=10k+7 \ .
$$
Moreover, $X_k$ is simply connected by the Lefschetz hyperplane theorem, and the second Betti number is given by $b_2(X_k)=c_2(X_k)-2=80k+74$.

On the other hand, consider the following family of Horikawa surfaces $Y_i$ from \cite{horikawa76}. 
Let $\Sigma_i$ be the Hirzebruch surface of degree $i$, that is the $\CP^1$-bundle over $\CP^1$ 
whose zero section $\Delta$ has self-intersection $-i$. Let $F$ denote the class of the fibre of the fibration.
Then we can construct the Horikawa surface $Y_i$ as the double cover $pr\colon Y_i\rightarrow\Sigma_i$ 
with branch locus homologous to $B=6\Delta+2(2i+3)F$.
Notice that these surfaces have ample canonical bundle $K_{Y_i}$ since 
$K_{Y_i}=pr^*(K_{\Sigma_i}+\frac{1}{2}B)=pr^*(\Delta+(i+1)F)$ and $\Delta+(i+1)F$ is an ample bundle. 
Moreover, $Y_i$ is simply connected because the branch locus $B$ is ample.
Now the characteristic numbers of $Y_i$ are:
$$
c_1^2(Y_i)=2i+4 \ ,\hspace{3mm}c_2(Y_i)=10i+56 \ ,\hspace{3mm}\chi(\mathcal{O}_{Y_i})=i+5 \ .
$$
Hence $b_2(Y_i)=10i+54$ so $b_2(X_k)=b_2(Y_i)$ for $i=8k+2$. 

From now on we denote $Y_{8k+2}$ by $Z_k$. For this we have:
$$
c_1^2(Z_k)=16k+8 \ ,\hspace{3mm}c_2(Z_k)=80k+76 \ ,\hspace{3mm}\chi(\mathcal{O}_{Z_k})=8k+7 \ .
$$ 
Both $X_k$ and $Z_k$ have ample canonical line bundle. 
Hence, by Proposition~\ref{SmaleGeiges}  the Boothby--Wang construction gives us two negative Sasaki structures with 
homotopic almost contact structures on the $(80k+73)$-fold connected sum $\#(80k+73)(S^2\times S^3)$.

Since the Sasaki structures are regular, their basic Hodge numbers are the Hodge numbers of the base surfaces.
Therefore we have
\begin{align*}
h^{0,2}(X_k)=10k+6,&\hspace{4mm}h^{1,1}(X_k)=60k+62,\\
h^{0,2}(Z_k)=8k+6,&\hspace{4mm}h^{1,1}(Z_k)=64k+62\ .
\end{align*}

Notice, on the one hand, that $\di(c_1(X_k))=\gcd\{k,2\}$. 
On the other hand, the main result of~\cite{nagami00} implies that $Z_k$ is spin if and only if $B/2$ is the 
Poincar\'e dual of the second Stiefel-Whitney class $w_2(\Sigma_{8k+2})$. In other words, $Z_k$ is spin 
if and only if $B/2$ is divisible by $2$ because $\Sigma_{8k+2}$ is spin. Now the intersection number of 
$B/2$ with $F$ equals $3$ and this implies that $Z_k$ is not spin.
In particular,  $\di(c_1(Z_k))$ is always odd, and so the two Sasaki structures do not have equivalent 
contact structures by Proposition~\ref{ThmHamilton} whenever $k$ is even.
\end{proof}

We now move on to the case of non-spin $5$-manifolds considered in Theorem~\ref{5nonspin}.

\begin{proof}[Proof of Theorem~\ref{5nonspin}]
Again we start with a $k$-tuple of surfaces $X_1,\ldots , X_k$ as in Theorem~\ref{surfaces}, 
so they are simply connected with ample canonical bundles, all with the same Euler characteristics,
but with pairwise distinct Hodge numbers. Moreover we now want all the $X_i$ to be non-spin. 
Inspecting the formula~\eqref{c1} for $c_1$ we see that it is enough to assume that 
$3-3q$ is always odd, for this ensures that $c_1$ is not divisible by $2$. Now the only constraint
we imposed on $q$ was that the different values of $3q-1$ should be pairwise coprime. We may
at the same time assume that all the $3q-1$ are odd, and so $3-3q$ is odd as well.

On $\C P^1\times \C P^2$ let $x_1$ and $x_2$ be the generators of the cohomology coming from 
the two factors. All classes of the form $e= ax_1+bx_2$ with positive $a$ and $b$ are K\"ahler 
classes. We fix such a class $e$ with $b$ even and $a$ and $b$ relatively prime. On every $X_i$ 
we perform the Boothby--Wang construction using as Euler classes the restrictions of $e$ to the $X_i$. 
This yields negative Sasaki structures on the total spaces. The Euler class is primitive because $a$ 
and $b$ are relatively prime, and so the total spaces are all simply connected with torsion-free
cohomology by Lemma~\ref{lem:spin}. Moreover, the assumption that all the $3q-1$ are odd but $b$ 
is even ensures $e$ and $c_1(X_i)$ have different reductions modulo $2$, and so the total spaces are 
not spin. Now Barden's classification~\cite{barden65} implies that all the total spaces are diffeomorphic 
to some $N_n=n(S^2\times S^3)\# (S^2\widetilde\times S^3)$, where $S^2\widetilde\times S^3$ is the non-trivial 
orientable $S^3$-bundle over $S^2$. Moreover, $n$ is the same for all our total spaces because all the 
$X_i$ have the same second Betti number. However, the fact that the $X_i$ have pairwise distinct Hodge 
numbers means that the basic Hodge numbers of the Sasaki structures on the total spaces are also 
pairwise distinct. This completes the proof.
\end{proof}

\section{Higher dimensions}

In this section we give examples of manifolds in dimensions $\geq 7$ admitting pairs of Sasaki
structures with different basic Hodge numbers. In particular we prove Theorem~\ref{main2}.

\begin{prop}\label{Thm5}
There exist closed $7$-manifolds admitting two negative Sasaki structures with different basic 
Hodge numbers. Moreover, manifolds with such pairs range over infinitely many homotopy types.
\end{prop}
\begin{proof}
Let $X_1$ and $X_2$ be two diffeomorphic $3$-dimensional complete intersections with different first Chern 
classes as given in Table~\ref{TableDiffeoHodgeNumbers}. After conjugating one of the complex structures if 
necessary we may assume that the diffeomorphism matches up the positive generators $x$ in 
$H^2(X_1)$ and $H^2(X_2)$. This generator is of course a \K class.
We perform the Boothby--Wang construction with Euler class $k \cdot x$ on both $X_i$.
The diffeomorphism between $X_1$ and $X_2$ then lifts to a diffeomorphism of the 
resulting negative Sasaki manifolds $(M_1,\eta_1,\phi_1,R_1,g_1)$ and $(M_2,\eta_2,\phi_2,R_2,g_2)$.
Lemma~\ref{LemHodgeChernCompInt} and Table~\ref{TableDiffeoHodgeNumbers} imply that these two
Sasaki structures are distinguished by their basic Hodge numbers which coincide with the Hodge numbers of $X_1$ and $X_2$.
As $k$ varies these examples range over infinitely many homotopy types since the fundamental group of $M_i$ is 
cyclic of order $k$ by the homotopy exact sequence of the Boothby--Wang fibration.
\end{proof}

For the remaining dimensions we give a uniform argument in the following theorem.
\begin{theor}\label{main3}
For all $n\geq 4$ there exist closed $(2n+1)$-dimensional manifolds 
admitting two negative Sasaki structures with different basic Hodge numbers.
Moreover, in each dimension manifolds with 
such pairs range over infinitely many homotopy types.
\end{theor}
\begin{proof}
Let $X_1$ and $X_2$ be two diffeomorphic $3$-dimensional complete intersections with different first Chern 
classes given in Table~\ref{TableDiffeoHodgeNumbers}.
Recall that $X_1$ and $X_2$ have different Hodge numbers by Lemma~\ref{LemHodgeChernCompInt}.

For any smooth projective algebraic variety $P$ of complex dimension $n-3$ with ample canonical bundle
the products $X_1\times P$ and $X_2\times P$ are diffeomorphic, and the diffeomorphism matches up the
positive primitive integral classes $x+K_{P}$, where $x$ is the positive generator of $H^2(X_i)$.
Therefore the
diffeomorphism lifts to the total spaces of the Boothby--Wang fibrations with these primitive classes
as Euler classes. These total spaces carry negative  Sasaki structures whose basic Hodge numbers
are the usual Hodge numbers of the bases $X_i\times P$. The compatibility of the Hodge and K\"unneth
decompositions shows that these Hodge numbers are different because those of the $X_i$ are.

Finally we can let $P$ range over infinitely many homotopy types, for example by using hypersurfaces 
of different degrees in $\C P^{n-2}$. This ensures that we obtain infinitely many homotopy types of 
Sasaki manifolds.
\end{proof}

We have now proved Theorem~\ref{main2}, since the case $n=2$ was covered in Theorems~\ref{main1}
and \ref{5nonspin}, the case $n=3$ was covered in Proposition~\ref{Thm5}, and the case $n\geq 4$
in Theorem~\ref{main3}. Moreover, the manifolds considered are simply connected unless $n=3$, where
 the fundamental groups are finite cyclic, and $n=4$ where the fundamental groups are those of 
 algebraic curves of general type.

\end{document}